%% file: jointSFS_arXiv.tex
\documentclass[aap,preprint,autosecdot]{imsart}

\RequirePackage[OT1]{fontenc}
\RequirePackage[numbers,sort]{natbib}
\RequirePackage[colorlinks,citecolor=blue,urlcolor=blue]{hyperref}   

\usepackage{graphicx}
\usepackage{amsmath,amssymb,amsthm,euscript,mathrsfs,bbm,caption,nicefrac}
\usepackage{subcaption}
\usepackage{mathtools,multirow}
\usepackage{array}
\usepackage{comment}
\usepackage[utf8]{inputenc} 
 \usepackage{appendix}
 \usepackage{cite}             
 \usepackage{import}
 \usepackage{nameref}

\newtheorem{proposition}{\noindent Proposition}
\newtheorem{lemma}{\noindent Lemma}

\def\rise#1#2{(#1)_{#2\uparrow}}
\def\I{\mathbb{I}}
\def\P{\mathbb{P}}
\def\E{\mathbb{E}}

\def\Coal{\mathcal{C}}
\def\CKill{\mathcal{K}}

\def\AncCoal{A^{\Coal}}
\def\AncCoaln#1{A^{\Coal^{#1}}}
\def\AncKill{A^{\CKill}}
\def\usfs{f}

\def\killedSet{K}
\def\killedSize{k}

\def\truncTime{\tau}
\def\timeVar{t}
\def\numAnc{m}

\def\mutSet{\mathcal{M}}

\def\popVertex{v}

\def\len{\mathcal{L}}
\def\numLoci{L}
\def\likelihood{\ell}

\def\tmrca{T_\text{\scriptsize MRCA}}

\def\bfmath#1{\mathchoice
        {\mbox{\boldmath$#1$}}%
        {\mbox{\boldmath$#1$}}%
        {\mbox{\boldmath$\scriptstyle#1$}}%
        {\mbox{\boldmath$\scriptscriptstyle#1$}}}%

\newcommand\truncWait[1]{T_{#1}^{\truncTime}}

\newcommand\polyaProb[4]{p_{#1,#2}^{#3,#4}}

\renewcommand\cite[1]{\citep{#1}}

\newcommand\lemref[1]{Lemma~\ref{#1}}

\newcommand\dt{\text{\rm d}t}
\newcommand\dx{\text{\rm d}x}

\begin{document}

\begin{frontmatter}
\title{Efficient computation of the joint sample frequency spectra for multiple populations}
\runtitle{}

\begin{aug}
\author{\fnms{John A.} \snm{Kamm}\thanksref{t1}\ead[label=e1]{jkamm@stat.berkeley.edu}},
\author{\fnms{Jonathan} \snm{Terhorst}\thanksref{t2}\ead[label=e2]{terhorst@stat.berkeley.edu}}
\and\\
\author{\fnms{Yun S.} \snm{Song}\thanksref{t1,t3}}
\ead[label=e3]{yss@stat.berkeley.edu}

\thankstext{t1}{Supported in part by an NIH grant R01-GM109454.}
\thankstext{t2}{Supported in part by a Citadel Fellowship.}
\thankstext{t3}{Supported in part by a Packard Fellowship for Science and Engineering, and a Miller Research Professorship.}

\runauthor{Kamm, Terhorst and Song}

\affiliation{University of California, Berkeley}

\address{J.~A. Kamm\\
Department of Statistics\\
University of California, Berkeley\\
Berkeley, CA 94720\\
USA\\
\printead{e1}}

\address{J. Terhorst\\
Department of Statistics\\
University of California, Berkeley\\
Berkeley, CA 94720\\
USA\\
\printead{e2}}

\address{Y.~S. Song\\
Department of Statistics and\\
\hspace{3mm} Computer Science Division\\
University of California, Berkeley\\
Berkeley, CA 94720\\
USA\\
\printead{e3}}  
\end{aug}

\begin{abstract}
  A wide range of studies in population genetics have employed the sample frequency spectrum (SFS), a summary statistic which   describes the distribution of mutant alleles at a polymorphic site in a sample of DNA sequences.
  In particular, recently there has been growing interest in analyzing the joint SFS data from multiple populations to infer parameters of complex demographic histories, including variable population sizes, population split times, migration rates, admixture proportions, and so on.
  Although much methodological progress has been made, existing SFS-based inference methods suffer from numerical instability and high computational complexity when multiple populations are involved and the sample size is large. 
  In this paper, we present new analytic formulas and algorithms that
  enable efficient computation of the expected joint SFS for multiple populations related by a complex demographic model with arbitrary population size histories (including piecewise exponential growth).
  Our results are implemented in a new software package called \emph{momi} (\underline{MO}ran \underline{M}odels for \underline{I}nference).
  Through an empirical study involving tens of populations, we demonstrate our improvements to numerical stability and computational complexity.
  \end{abstract}



\end{frontmatter}

\section{Introduction}

The sample frequency spectrum (SFS) is the distribution of allele frequencies at a polymorphic site in a collection of DNA sequences randomly drawn from a population.  This summary statistic is used in a variety of inference problems in population genetics \citep{wakeley1997estimating,griffiths1998age,nielsen2000estimation,gutenkunst2009inferring,coventry2010deep,gazave2014neutral,gravel2011demographic, 
nelson2012abundance,excoffier2013robust,jenkins2014general,bhaskar2015efficient}, often in the context of likelihood-based analysis of single nucleotide polymorphism (SNP) data.
Over the past several years, there has been much interest in analyzing the joint SFS data from multiple populations to infer complex demographic models involving population size changes, population splits, migration, and admixture.
Inferring population demographic histories is not only intrinsically interesting, for example in dating events such as the out-of-Africa migration of modern humans \citep{schaffner2005calibrating, gutenkunst2009inferring}, but
is also important for biological applications, such as distinguishing between the effects of natural selection and
demography \citep{beaumont1996evaluating,boyko2008assessing}.

Likelihood-based inference methods using the SFS require accurate computation of the expected SFS under a given demographic model.
As further detailed below, however, existing methods suffer from numerical instability and high computational complexity when multiple populations are involved and the sample size is large. 
The joint SFS for multiple populations describes the distribution of joint allele frequencies across the different populations.
In this paper, we tackle the problem of computing the expected joint SFS for many populations, given a complex demographic model relating them.

The SFS has been studied in the context of two dual processes, the
Wright-Fisher diffusion \citep{kimura1955solution} and Kingman's coalescent
\citep{fu1995statistical}, and both approaches can be used to compute the multi-population
SFS. In the diffusion approach of \citet{gutenkunst2009inferring}, which was later further extended \citep{gravel2011demographic,
lukic2012demographic}, one numerically solves partial differential equations forward in time to approximate the distribution of joint allele frequencies at present. The diffusion framework has the advantage of being applicable to arbitrary demographic models, but its computational complexity grows
exponentially with the number of populations, and current implementations have
difficulty handling more than three \citep{gutenkunst2009inferring} or four
populations \citep{lukic2012demographic}.

In the coalescent approach, the SFS is computed by integrating over all
genealogies underlying the sample.  This can be done via Monte
Carlo or analytically.  Monte Carlo integration approach \citep{nielsen2000estimation} can effectively handle arbitrary demographic
histories with a large number of populations, and 
\citet{excoffier2013robust} have recently developed a useful implementation.
However, when the number $\mathcal{D}$ of populations (or demes) is
moderate to large, most of the $O(n^\mathcal{D})$ SFS entries, where $n$ denotes the sample size, will be unobserved in
simulations, and thus the Monte Carlo integral may naively assign a probability
of $0$ to observed SNPs. Monte Carlo computation of the SFS thus requires
careful regularization techniques to avoid degeneracy issues.

An alternative to  the Monte Carlo approach is to compute the
SFS exactly via analytic integration over coalescent genealogies \citep{wakeley1997estimating,griffiths1998age}.
For a demography involving multiple populations, this can be done efficiently by a dynamic program \citep{chen2012joint, chen2013intercoalescence}.
This algorithm is more complicated and less general than both the Monte
Carlo and diffusion approaches: while it can handle population splits, merges,
size changes, and instantaneous gene flow, it is difficult to include continuous gene flow between populations.
However, it scales well to a large number $\mathcal{D}$ of populations, since it only computes entries of the SFS that are observed in
the data, and ignores the $O(n^\mathcal{D})$ SFS entries that are not
observed.
Unfortunately, existing coalescent-based algorithms \citep{wakeley1997estimating,chen2012joint, chen2013intercoalescence} do not scale well to a large sample size $n$,
both in terms of running time and numerical stability.
In particular, the algorithm relies on large alternating sums that
explode with $n$ and exhibit catastrophic cancellation.

In this paper, we significantly improve the computational complexity
and numerical stability of the coalescent approach.
We show how the alternating sums can be avoided altogether, and replaced with faster and more stable formulas.
Moreover, we introduce a second speedup by replacing the coalescent with a Moran model in the dynamic program.

The dynamic program algorithm to compute the SFS involves splitting the demography
into its component subpopulations, each of which contains
a single population coalescent, but truncated at some time
$\truncTime$ in the past.
In Section~\ref{sec:truncSfs}, we focus on this \emph{truncated coalescent}.
In particular, we focus on computing the \emph{truncated SFS} $\usfs_n^\truncTime(\killedSize)$, the expected number of mutations arising in the time interval $[0,\truncTime)$ which subtend exactly $k$ out of $n$ individuals sampled at time $0$.
We give an algorithm for computing $\usfs_n^\truncTime(\killedSize)$ efficiently,
using recurrence relations combined with results from \citet{polanski2003note,polanski2003new} and \citet{bhaskar2015efficient}.
We also provide an alternative formula for $\usfs_n^\truncTime(\killedSize)$ based on  \emph{the coalescent with killing.}

In Section~\ref{sec:jointSfs}, we describe the
coalescent algorithm of \citet{chen2012joint, chen2013intercoalescence},
and show how our formulas for $\usfs_n^\truncTime(\killedSize)$
improve its computational complexity.
For the special case where the demographic history forms a tree,
we introduce an additional speedup by replacing the coalescent
with a Moran model.
For such tree-shaped demographies,
we can compute the observed SFS entries in $O(n^2
  \mathcal{D} + n \log(n) \mathcal{D} \numLoci)$, where $n$ is the sample size,
  $\mathcal{D}$ is the number of populations at the present, and $\numLoci$ is the
  number of observed entries in the SFS. This is an improvement over the $O(n^5
  \mathcal{D} + n^4 \mathcal{D} \numLoci)$ complexity in \citet{chen2012joint,
  chen2013intercoalescence}.
For more general demographic histories with migration or admixture,
the algorithm of \citet{chen2012joint, chen2013intercoalescence} is
$O(n^5 V + W \numLoci)$,
where $V$ is the number populations (vertices) throughout the history,
and $W$ is a term that depends on $n$ and the graph structure of the
demography;
we improve this to $O(n^2 V + W \numLoci)$.
In future work, we will give explicit expressions for $W$,
and extend our Moran-based speedup to demographies with pulse migration.

We note that some of our improvements
are related to results in \citet{bryant2012inferring},
whose $O(n^2 \log(n) \mathcal{D} \numLoci)$ algorithm
computes the one-locus likelihood for 
species trees with recurrent mutation and piecewise
constant population sizes.
By contrast, our method, like that of
\citet{chen2012joint, chen2013intercoalescence},
considers an infinite sites model \citep{kimura1969number} without recurrent
mutation, but can handle exponentially growing population sizes.
In fact, our method goes even further, and easily accommodates arbitrary population size changes.

In Section~\ref{sec:results}, we demonstrate
the improved speed and accuracy of our algorithm in a numerical study involving tens of populations.
We implement and release our algorithm in a publicly available Python package,
\emph{momi} (\underline{MO}ran \underline{M}odels for \underline{I}nference).
Proofs of the mathematical results presented in Section~\ref{sec:truncSfs} are provided in Section~\ref{sec:proofs},

\section{The truncated sample frequency spectrum}
\label{sec:truncSfs}

\subsection{Background}
We denote Kingman's coalescent~\citep{kingman1982coalescent,kingman1982genealogy,kingman1982exchangeability} on $n$ leaves
$\{\Coal^n_\timeVar\}_{\timeVar \geq 0}$ to be the backward-in-time
Markov jump process, whose
value at time $\timeVar$ is a partition of
$\{1,\ldots,n\}$,
and at time $\timeVar$, each pairs of blocks in $\Coal^n_\timeVar$
coalesce with rate $\alpha(\timeVar)$.
We also call $\frac1{\alpha(\timeVar)}$ the \emph{population size history function}.
We denote the ancestral process $\AncCoaln{n}_\timeVar = |\Coal^n_\timeVar|$
to be the number of blocks in $\Coal^n_\timeVar$,
so that $\AncCoaln{n}_\timeVar$ is a pure death process
with $\AncCoaln{n}_0 = n$ and the rate of transition from $\numAnc$ to $\numAnc-1$ given by
$\lambda^{\Coal}_{\numAnc,\numAnc-1}(\timeVar) = {\numAnc \choose 2} \alpha(\timeVar)$.

We often drop the dependence
on $n$,
and write $\Coal_\timeVar = \Coal^n_\timeVar$
and $\AncCoal_\timeVar = \AncCoaln{n}_\timeVar$.
We prefer to denote a dependence on
$n$ through the
probability $\P_n$
and the expectation $\E_n$.
So if $X(\Coal^n)$ denotes
a random variable of the process
$\Coal^n$,
we usually write $\E_n[X]$
instead of $\E[X(\Coal^n)]$.

Let $\xi_i$ denote the partition
of $\{1,\ldots,n\}$ when $\Coal_\timeVar$
has $i$ blocks (also referred to as lineages).
Let $T_i = \int_0^\infty \I_{\AncCoal_\timeVar = i} d\timeVar$
denote the amount of time $\Coal_\timeVar$
has exactly $i$ lineages.
It is a fundamental fact of the coalescent
that the waiting times $\mathbf{T}_{n:2} = (T_n, \ldots, T_2)$
are independent of the partitions $\bfmath{\xi}_{n:2} = (\xi_n,
\ldots, \xi_2)$ \citep{kingman1982coalescent}.

A sample path of $\Coal^n$
can be viewed as a rooted ultrametric tree with $n$ leaves
labeled $1,\ldots,n$,
where $\Coal_\timeVar$ is the partition induced
on $\{1,\ldots,n\}$
by cutting the tree at height $\timeVar$.
Now suppose we drop mutations
onto this tree as a Poisson point process
with rate $\frac{\theta}2$,
and let
$\mutSet$ denote the set of leaves
that are beneath mutations (where we only
consider mutations beneath the root, so by assumption $\mutSet \neq
\{1,\ldots,n\}$).
Then we define the sample frequency spectrum
$\usfs_n(\killedSize)$, for $0 < \killedSize < n$, as
the first order Taylor series coefficient of $\P_n(|\mutSet| =
\killedSize)$
in the mutation rate,
\begin{eqnarray*}
  \P_n(|\mutSet| = \killedSize) &=& \frac{\theta}2 \usfs_n(\killedSize)
  + o(\theta).
\end{eqnarray*}
We will generally refer to $\usfs_n(\killedSize)$ as the sample frequency
spectrum (SFS).

We also note two alternative definitions of the SFS.
First, $\usfs_n(\killedSize)$ is the expected number
of mutations with $\killedSize$ descendants
when $\frac{\theta}2 = 1$.
Second, $\frac{1}{{n \choose  |\killedSet|}}\usfs_n(|\killedSet|)$ is the expected length of the branch
whose leaf set is $\killedSet$.
More specifically, let $\I$ denote the indicator function,
and define $\len_\killedSet := \int_0^\infty \I_{\killedSet \in \Coal_\timeVar}
d\timeVar$.
Then
\begin{eqnarray*}
\frac{1}{{n \choose  |\killedSet|}}\usfs_n(|\killedSet|) &=& \E_n[\len_\killedSet].
\end{eqnarray*}
The equivalence of these alternate definitions
follows from
 previous results in
\citet{griffiths1998age, jenkins2011effect, bhaskar2012approximate}.

Note the SFS is sometimes defined to be a normalized version
of $\usfs_n(\killedSize)$, so that the entries sum to 1. We do not
follow that convention, and use the unnormalized definition
for the SFS throughout this paper.

\begin{figure}[t]
  \centering
  \def\svgwidth{.5\textwidth}
  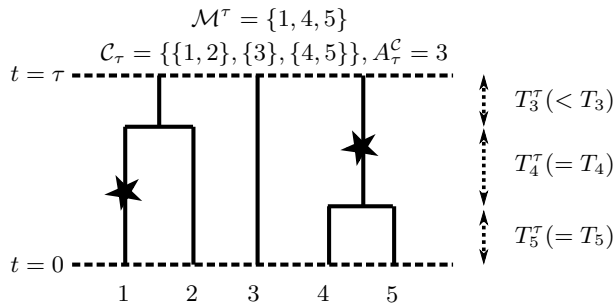
  \caption{A sample path of the coalescent truncated at time $\truncTime$.  
  Star symbols denote mutations, while $\mutSet^\truncTime$ denotes the set of leaves under those mutations.
  $T_k^\truncTime$ denotes the waiting time in the interval $[0,\truncTime)$ while there are $k$ lineages.
  }
  \label{fig:notation}
\end{figure}

\subsection{The truncated coalescent and SFS}

We now consider truncating the coalescent with mutation
at time $\truncTime$, as illustrated in Figure~\ref{fig:notation}.
Let $\mutSet^\truncTime$ denote
the set of leaves under mutations
occurring in the time interval $[0,\truncTime)$.
We define the \emph{truncated} SFS
$\usfs_n^\truncTime(\killedSize)$ according to
\begin{eqnarray*}
  \P_n(|\mutSet^\truncTime| = \killedSize) &=& \frac{\theta}2
  \usfs_n^\truncTime(\killedSize) + o(\theta).
\end{eqnarray*}
By the same arguments
as for the untruncated SFS, one can show that
$\usfs_n^\truncTime(\killedSize)$ gives
the expected number of mutations in $[0,\truncTime)$
with $\killedSize$ descendants,
and letting
$\len^\truncTime_\killedSet := \int_0^\truncTime \I_{\killedSet
  \in \Coal_\timeVar} d\timeVar$
denote the branch length subtending $\killedSet\subset\{1,\ldots,n\}$ within
$[0,\truncTime)$,
we have
\begin{eqnarray*}
  \frac{1}{{n \choose |\killedSet|}}\usfs_n^\truncTime(|\killedSet|) &=& \E_n[\len_\killedSet^\truncTime].
\end{eqnarray*}
Note that for $\killedSize < n$,
we have
$\usfs_n(\killedSize) = \usfs_n^\infty(\killedSize)$.
For the truncated SFS, we will also consider
mutations above the root,
and so allow $k=n$ (i.e., $\mutSet^\truncTime = \{1,\ldots,n\}$),
with $\usfs_n^\truncTime(n) = \E_n[\len_{\{1,\ldots,n\}}^\truncTime]$ giving the expected number of
mutations within $[0,\truncTime)$ subtending the whole sample.

Given a random variable $X$,
we define conditional versions of
the SFS $\usfs_n^\truncTime(\killedSize \mid X)$
according to
\begin{eqnarray*}
  \P_n(|\mutSet^\truncTime| = \killedSize \mid X) &=& \frac{\theta}2
  \usfs_n^\truncTime(\killedSize \mid X) + o(\theta).  
\end{eqnarray*}

An example of a useful conditional SFS
is $\usfs_n^\truncTime(\killedSize \mid \AncCoal_\truncTime =
\numAnc)$,
the expected branch length subtending
$\killedSize$ leaves
given $\numAnc$ ancestors at time $\truncTime$.
In particular, \citet{chen2012joint}
devised a dynamic program algorithm
to compute the joint SFS for multiple populations
under complex demographic histories,
by computing
$\{\usfs_\nu^\truncTime(\killedSize)\}_{\killedSize \leq \nu \leq n}$
on each subpopulation of the history,
where $\truncTime$ is the length of time a particular subpopulation
exists. The unconditional SFS
$\usfs_\nu^\truncTime(\killedSize)$
is in turn computed
in terms of
$\usfs_\nu^\truncTime(\killedSize \mid \AncCoal_\truncTime = \numAnc)$
by writing
\begin{eqnarray}
  \usfs_\nu^\truncTime(\killedSize)
&=&
\sum_{m=1}^{n-\killedSize+1} 
\P_\nu(\AncCoal_\truncTime = \numAnc)
\usfs_\nu^\truncTime(\killedSize \mid \AncCoal_\truncTime =
\numAnc) . \label{eq:tsfs}
\end{eqnarray}
In Section~\ref{sec:complex:sfs}, we describe
the dynamic program algorithm for
computing
the joint SFS for multiple populations,
and the way in which this algorithm uses
the terms $\usfs_\nu^\truncTime(\killedSize)$.

We consider how to compute
\eqref{eq:tsfs}.
The first term in the summand, 
$\P_\nu(\AncCoal_\truncTime = \numAnc)$,
can be computed in at least three ways:
by numerically exponentiating
the rate matrix of $\AncCoal$,
by computing an alternating sum with $O(\nu)$ terms
\citep{tavare1984line},
or by solving a recursion we describe
in Section~\ref{sec:ancProbs}.
We note that the recursion
described in Section~\ref{sec:ancProbs}
has the 
advantage of computing
all values of
$\P_\nu(\AncCoal_\truncTime = \numAnc)$,
$\numAnc \leq \nu \leq n$,
in $O(n^2)$ time.

The second term 
$\usfs_\nu^\truncTime(\killedSize \mid \AncCoal_\truncTime =
\numAnc)$ in the summand of \eqref{eq:tsfs} is computed in 
\citet{chen2012joint}
as
\begin{align}
  \usfs_\nu^\truncTime(\killedSize \mid \AncCoal_\truncTime = \numAnc) &=
  \sum_{i=\numAnc}^{\nu}  i \polyaProb{\nu}{i}{\killedSize}{1}
  \E_\nu[\truncWait{i} \mid \AncCoal_\truncTime =
  \numAnc ],\label{eq:csfs:numAnc}
\end{align}
where
\begin{align*}
  \polyaProb{\nu}{i}{\killedSize}{j} &:= 
  \begin{cases}
    \frac{{\killedSize - 1 \choose j - 1}{\nu - \killedSize - 1 \choose
        i - j - 1}}{{\nu - 1 \choose i - 1}}, & \text{if $\killedSize
      \geq j > 0$ and $\nu - \killedSize \geq i -
    j > 0$}, \\
    1, & \text{if $j = \killedSize = 0$ or $i - j = \nu - \killedSize =
      0$}, \\
    0, & \text{else},
  \end{cases}
\end{align*}
is the transition probability of the P\'olya urn model,
starting with $i - j$ white balls and $j$ black balls,
and ending with $\nu - \killedSize$ white balls and $\killedSize$ black
balls \citep{johnson1977urn},
and
\begin{align*}
\truncWait{i} := \int_0^\truncTime \I_{\AncCoal_\timeVar = i} d\timeVar  
\end{align*}
is the length of time in $[0,\truncTime)$ where there are $i$
ancestral lineages to the sample, as illustrated in Figure~\ref{fig:notation}.
\citet{chen2012joint} provides a formula
for the conditional expectation
$\E_\nu[\truncWait{i} \mid \AncCoal_\truncTime =
\numAnc]$
for the case of constant population size, which he later extends \citep{chen2013intercoalescence} to the case of an exponentially growing
population.
However, these formulas involve
a large alternating sum with $O(\nu^2)$ terms.
Thus, computing
$\E_\nu[\truncWait{i} \mid \AncCoal_\truncTime =
\numAnc]$
for every value of $i,m,\nu$,
as required to compute
$\{\usfs_\nu^\truncTime(\killedSize)\}_{\killedSize \leq \nu \leq n}$
with
\eqref{eq:tsfs} and \eqref{eq:csfs:numAnc},
takes $O(n^5)$ time
with these formulas.
In addition, large alternating sums
are numerically unstable due to catastrophic
cancellation
\citep{higham2002accuracy},
and so these formulas
require the use of high-precision
numerical libraries, further
increasing runtime.

\subsection{An efficient, numerically stable algorithm for computing the truncated SFS}

Here, we present a numerically stable algorithm to compute, for a given positive integer $n$,  all of $\{\usfs_\nu^\truncTime(\killedSize) \mid 1\leq k\leq \nu\leq n\}$ in $O(n^2)$ time instead of $O(n^5)$ time.
Our approach utilizes the following two lemmas:
\begin{lemma}
\label{thm:tsfs:tmrca}
The entry $\usfs_n^\truncTime(n)$ of the truncated SFS is
given by
\begin{eqnarray}
\usfs_n^\truncTime(n)
&=&
\tau - \sum_{\killedSize=1}^{n-1} \frac{\killedSize}n \usfs_n^\truncTime(\killedSize).
\label{eq:tsfs:tmrca}
\end{eqnarray}
\end{lemma}
\begin{lemma}
\label{thm:tsfs:recurrence}
For all $1\leq k\leq \nu$, the truncated SFS $\usfs_\nu^\truncTime(\killedSize)$ satisfies
the linear recurrence
\begin{eqnarray}
\usfs_\nu^\truncTime(\killedSize)
&=&
\frac{\nu-k+1}{\nu+1} \usfs_{\nu+1}^\truncTime(\killedSize)
+
\frac{k+1}{\nu+1} \usfs_{\nu+1}^\truncTime(\killedSize+1).
\label{eq:tsfs:recurrence}
\end{eqnarray}
\end{lemma}
We prove \lemref{thm:tsfs:tmrca} in Section~\ref{proof:tsfs:tmrca}.
We note here that our proof also yields the identity $\E[\tmrca] =
\sum_{\killedSize=1}^{n-1} \frac{\killedSize}n \usfs_n(\killedSize)$, where 
$\tmrca$ denotes the time to the most recent common ancestor of the
sample; to our knowledge, this formula is new.
A proof of \lemref{thm:tsfs:recurrence} is provided in Section~\ref{proof:tsfs:recurrence}.

We now sketch our algorithm.   For a given $n$, we show below that all
values of $\usfs_n^\truncTime(\killedSize)$, for $1 \leq k < n$,
can be computed in $O(n^2)$ time.
We then compute $\usfs_n^\truncTime(n)$ using
\lemref{thm:tsfs:tmrca} in $O(n)$ time.
Finally, using $\usfs_n^\truncTime(\killedSize)$ for $1 \leq
\killedSize \leq n$ as boundary conditions, 
\lemref{thm:tsfs:recurrence} allows us to compute all $\usfs_\nu^\truncTime(\killedSize)$, for $\nu =  n-1,n-2,\ldots,2$ and $k = 1, \ldots, \nu$, in $O(n^2)$ time.

We now describe how to compute the aforementioned terms
$\usfs_n^\truncTime(\killedSize)$, for all $\killedSize < n$, in $O(n^2)$ time.
We first recall the result of \citet{polanski2003new} which represents the untruncated SFS $\usfs_n(\killedSize)$, for $1 \leq k \leq n-1$, as 
\begin{equation}
  \usfs_n(\killedSize) = \sum_{m=2}^n W_{n,k,m} c_m, 
\label{eq:PK}
\end{equation}
where 
\begin{align}
c_m := \E_m[T_{m}] & = 
\int_0^\infty t {m\choose 2} \alpha(t) \exp\Bigg[-{m\choose 2} \int_0^t \alpha(x) \dx\Bigg] \dt \nonumber \\
& = \int_{0}^\infty \exp\Bigg[-{m\choose 2} \int_0^t \alpha(x) \dx\Bigg] \dt
\label{eq:cm}
\end{align}
denotes the waiting time to the first coalescence for a sample of size $m$, and $W_{n,k,m}$ are universal constants that are efficiently computable using the following recursions \citep{polanski2003new}:
\begin{align}
  W_{n,k,2} &= \frac{6}{n+1}, \nonumber \\
  W_{n,k,3} &= 30 \frac{(n - 2k)}{(n+1)(n+2)}, \nonumber\\
  W_{n,k,m+2} &= - \frac{(1 + m)(3 + 2m)(n - m)}{m(2m - 1)(n + m + 1)} W_{n,k,m} + \frac{(3 + 2m)(n - 2k)}{m(n + m + 1)} W_{n,k,m+1},\nonumber\\ 
  \label{eq:W_nkm}
\end{align}
for $2 \leq m \leq n - 2$.  The key observation is to note that, in a similar vein as \eqref{eq:PK}, we have:
\begin{lemma}
  \label{thm:truncated:polanski:kimmel}
The truncated SFS $\usfs_n^\truncTime(\killedSize)$, for $1 \leq k \leq n-1$, can be written as 
\begin{equation}
  \usfs_n^\truncTime(\killedSize) = \sum_{m=2}^n W_{n,k,m} c_m^\truncTime, 
  \label{eq:PK_truncated}
\end{equation}
where $c_m^\truncTime$ is a truncated version of \eqref{eq:cm}:
\begin{equation}
  c_m^\truncTime := \E_m[T_{m}^\truncTime] = \int_{0}^\truncTime  \exp\Bigg[-{m\choose 2} \int_0^t \alpha(x) \dx\Bigg] \dt.
  \label{eq:cm_truncated}
\end{equation}
\end{lemma}
We prove \lemref{thm:truncated:polanski:kimmel} in Section~\ref{proof:truncated:polanski:kimmel}.
For piecewise-exponential $\alpha(\timeVar)$, $c_m^\truncTime$ can be computed explicitly using formulas from \citet{bhaskar2015efficient}.
Using \eqref{eq:W_nkm}, we can compute all values of $W_{n,k,m}$, for $1\leq k \leq n$ and $2 \leq m \leq n$, in $O(n^2)$ time.  Then, using \eqref{eq:PK_truncated}, all values of $\usfs_n^\truncTime(\killedSize)$, for $1\leq k \leq n-1$ can be computed in $O(n^2)$ time.

Note that the above algorithm not only significantly improves computational complexity, but also resolves numerical issues, since it allows us to avoid computing the expected times $\E_\nu[\truncWait{i} \mid \AncCoal_\truncTime = \numAnc]$,
which are alternating sums of $O(n^2)$ terms and are numerically unstable to evaluate for large values of $n$ (say, $n>50$).

\subsection{An alternative formula for piecewise-constant subpopulation sizes}

For demographic scenarios with piecewise-constant subpopulation sizes, we present an alternative formula for computing the truncated SFS within a constant piece. This formula 
has the same sample computational complexity as that described in the previous section.

\begin{figure}[t]
  \centering
  \def\svgwidth{.5\textwidth}
  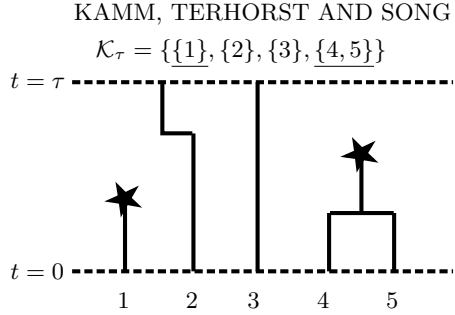
  \caption{The coalescent with killing for the genealogy
  in Figure~\ref{fig:notation}. Note that $\CKill_\truncTime$ is a
  marked partition, with the blocks killed by mutations in $[0,\truncTime)$ being
  specially marked.}
  \label{fig:coalkill}
\end{figure}

Let
$\CKill_\timeVar$
denote the \emph{coalescent with killing},
a stochastic process
that is closely related to
the Chinese restaurant process,
Hoppe's urn,
and Ewens' sampling formula
\citep{aldous1985exchangeability, hoppe1984polya}.
In particular,
the coalescent with killing $\{\CKill_\timeVar\}_{\timeVar \geq 0}$
is a stochastic process whose value at time $\timeVar$
is a \emph{marked} partition of $\{1,\ldots,n\}$,
where each partition block is marked as ``killed'' or
``unkilled''.
We obtain the partition for $\CKill_\timeVar$
by dropping mutations onto
the coalescent tree as a Poisson point process
with rate $\frac{\theta}2$,
and then defining
an equivalence relation
on $\{1,\ldots,n\}$,
where $i \sim j$
if and only if $i,j$ have coalesced by time $\timeVar$
and there are no mutations on the branches between $i$ and $j$
(i.e., $i$ and $j$ are identical by descent).
We furthermore mark the equivalence classes (i.e. partition blocks)
of $\CKill_\timeVar$
that are descended from a mutation in $[0,\timeVar)$
as ``killed''.
See Figure~\ref{fig:coalkill} for an illustration.
The process $\CKill_\truncTime$
can also be obtained by running
Hoppe's urn, or equivalently the
Chinese restaurant process,
forward in time
\citep[Theorem 1.9]{durrett2008probability}.

Let $\AncKill_\timeVar$
be the number of unkilled
blocks in $\CKill_\timeVar$,
so that $\AncKill_\timeVar$ is a pure death
process with transition rate $\lambda^\CKill_{i,i-1}(\timeVar) = {i \choose 2}
\alpha(\timeVar) + \frac{i \theta}2$
(the rate of coalescence is the number of unkilled pairs ${i
  \choose 2} \alpha(\timeVar)$,
and the rate of killing due to mutation is $\frac{i \theta}2$).
Our next proposition
gives a formula for the
truncated conditional sample frequency spectrum given
$\AncKill_\truncTime$, i.e., $\usfs_n^\truncTime(\killedSize \mid
\AncKill_\truncTime = \numAnc)$.

\begin{proposition} \label{thm:csfs:ancKill}
Consider the constant population size history $\frac1{\alpha(\timeVar)} = \frac1\alpha$
for $\timeVar \in [0,\truncTime)$, and let $\numAnc > 0$ and $0 < \killedSize \leq n - \numAnc$.
The joint probability that the number of derived mutants
is $\killedSize$
and  the number of unkilled ancestral lineages is
$\numAnc$, when truncating at time $\truncTime$, is
given by
  \begin{eqnarray*}
    \P_n(|\mutSet^\truncTime| = \killedSize, \AncKill_\truncTime =
    \numAnc)
&=&
  \frac{\theta}2 \usfs_n^\truncTime(\killedSize
  \mid \AncKill_\truncTime = \numAnc) \P(\AncCoal_\truncTime =
  \numAnc)  + o(\theta) ,
  \end{eqnarray*}
  where
  \begin{eqnarray}
    \usfs_n^\truncTime(\killedSize \mid \AncKill_\truncTime = \numAnc)
    &=&
    \frac2{\alpha\killedSize} \frac{{n - \numAnc \choose \killedSize}}{{n - 1
        \choose \killedSize}}.
    \label{eq:csfs:ancKill}
  \end{eqnarray}
\end{proposition}

We prove Proposition~\ref{thm:csfs:ancKill}
in Section~\ref{proof:csfs:ancKill}.
Note that this equation does not hold for the case
$\killedSize = n, \numAnc = 0$,
but fortunately we do not need to consider that case in what follows below.

We can use
Proposition~\ref{thm:csfs:ancKill}
to
stably and efficiently
compute the terms $\usfs_\nu^\truncTime(\killedSize)$,
for
$\killedSize \leq \nu \leq n$,
as follows.
We first compute the case
$\killedSize < \nu = n$.
Note that $\P_n(|\mutSet^\truncTime| = \killedSet) = \sum_\numAnc
\P_n(|\mutSet^\truncTime| = \killedSet, \AncKill_\truncTime = \numAnc)$.
So for $\killedSize < n$,
by Proposition~\ref{thm:csfs:ancKill}
\begin{eqnarray}
  \usfs_n^\truncTime(\killedSize)
&=&
\sum_{\numAnc=1}^n \usfs_n^\truncTime(\killedSize \mid
\AncKill_\truncTime = \numAnc)
\P_n(\AncCoal_\truncTime = \numAnc)
\notag
\\
&=&
\sum_{\numAnc=1}^n
    \frac2{\alpha\killedSize} \frac{{n - \numAnc \choose \killedSize}}{{n - 1
        \choose \killedSize}}
\P_n(\AncCoal_\truncTime = \numAnc). \label{eq:ancKilln}
\end{eqnarray}
The sum in \eqref{eq:ancKilln} contains $O(n)$ terms,
so it costs $O(n^2)$ to compute $\usfs_n^\truncTime(\killedSize)$
for all $\killedSize < n$.
After this, we use \lemref{thm:tsfs:tmrca} to compute $\usfs_n^\truncTime(n)$,
and then use \lemref{thm:tsfs:recurrence} 
to compute $\usfs_\nu^\truncTime(\killedSize)$ for all $1\leq \killedSize \leq \nu < n$.
Since there are $O(n^2)$ such terms, this also takes $O(n^2)$ time.

\section{The joint SFS for multiple populations}\label{sec:jointSfs}

In this section we discuss an algorithm
for computing the multi-population SFS \citep{wakeley1997estimating,chen2012joint,chen2013intercoalescence}.
We describe the algorithm in Section~\ref{sec:complex:sfs},
and note how the results from Section~\ref{sec:truncSfs}
improve the time complexity of this algorithm.
In Section~\ref{sec:moran}, we focus on the special case
of tree-shaped demographies, and introduce a further algorithmic
speedup by replacing the coalescent with a Moran model.

Let $V$ be the number of subpopulations in the demographic history,
$n$ the total sample size, and $L$ the number of SFS entries to
compute. Then the results from Section~\ref{sec:truncSfs} improve the computational complexity of the
SFS from $O(n^5 V + WL)$ to $O(n^2V + WL)$, where $W$ is a term that
depends on the structure of the demographic history.
In the special case of tree-shaped demographies, the algorithm from
\citet{chen2012joint} gives $W = O(n^4V)$. The Moran-based speedup
from Section~\ref{sec:moran}, combined with results from \citet{bryant2012inferring}, improves this to $W = O(n \log(n) V)$.

The Moran-based speedup can be generalized to non-tree demographies,
but the notation, implementation, and analysis of computational
complexity becomes substantially more complicated. We thus
leave its generalization to future work.

\begin{figure}[t]
  \centering
  \def\svgwidth{0.95\textwidth}
  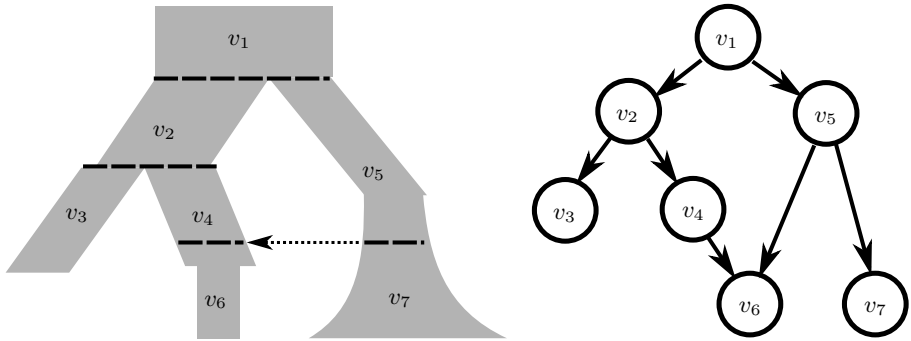
  \caption{A demographic history with a pulse migration event (left), and its corresponding
    directed graph (right).}
  \label{fig:demography}
\end{figure}

\subsection{A coalescent-based dynamic program\label{sec:complex:sfs}}

Suppose at the present we have
$\mathcal{D}$ populations,
and in the $i$th population
we observe
$n_i$ alleles.
For a single point mutation,
let $\mathbf{x} = (x_1, \ldots, x_\mathcal{D})$
denote the number of alleles
that are derived in each population.
We wish
to compute $\usfs(\mathbf{x})$,
where $\frac{\theta}2 \usfs(\mathbf{x})$
is
the expected number of
point mutations with derived
counts $\mathbf{x}$.

For demographic
histories
consisting
of population size changes,
population splits,
population mergers,
and pulse admixture events,
\citet{chen2012joint}
gave an algorithm to compute
$\usfs(\mathbf{x})$
using the truncated SFS
$\usfs_n^\truncTime(\killedSize)$
that we defined in Section~\ref{sec:truncSfs}.

We describe this algorithm 
to compute $\usfs(\mathbf{x})$.
We start
by representing the population history
as a directed acyclic graph (DAG),
where each vertex $\popVertex$ represents
a subpopulation (Figure~\ref{fig:demography}).
We draw a directed edge
from $\popVertex$ to $\popVertex'$
if there is gene flow
from the bottom-most part
of $\popVertex$
to the top-most part
of $\popVertex'$,
where ``down'' is the present
and ``up'' is the ancient past.
Thus, the leaf vertices
correspond to the subpopulations
at the present.
For a vertex $\popVertex$
in the population history graph,
let
$\truncTime_\popVertex \in (0,\infty)$
denote the length of time
the corresponding population persists,
and let $\alpha_\popVertex: [0,\truncTime_\popVertex) \to \mathbb{R}^+$
denote the inverse population size history of $\popVertex$.
So going backwards in time from the present,
$\alpha_\popVertex(\timeVar)$ gives the instantaneous rate at which
two lineages in $\popVertex$ coalesce,
after $\popVertex$ has existed for time $\timeVar$.
We use $\usfs^\popVertex_n(\killedSize)$ to
denote the truncated SFS for the coalescent
embedded in $\popVertex$,
i.e.,
$\usfs^\popVertex_n(\killedSize) = \usfs^{\truncTime_\popVertex}_n(\killedSize)$
for a coalescent with 
coalescence rate $\alpha_\popVertex(\timeVar)$.
Then we have
\begin{eqnarray}
  \usfs(\mathbf{x})
&=&
\sum_{\popVertex}
\sum_{\numAnc_0^\popVertex, \killedSize_0^\popVertex}
\usfs^\popVertex_{\numAnc_0^\popVertex}(\killedSize_0^\popVertex)
\P(\mathbf{x} \mid \killedSize_0^\popVertex, \numAnc_0^\popVertex)
\P(\numAnc_0^\popVertex)
\label{eq:chen:sfs:formula}
\end{eqnarray}
where $\numAnc_0^\popVertex$
denotes
the number of lineages 
at
the bottom of $\popVertex$
that are
ancestral to the initial sample,
and
$\killedSize_0^\popVertex$
denotes the number
of these lineages with a derived
allele.

In order to use
\eqref{eq:chen:sfs:formula},
we must compute
$\usfs^\popVertex_{\numAnc_0^\popVertex}(\killedSize_0^\popVertex)$
for every population $\popVertex$,
and every value of $\numAnc_0^\popVertex$
and $\killedSize_0^\popVertex$.
If $n$ is the total sample size
and $V$ the total number of vertices,
then this takes $O(n^5V)$ time
using the formulas of 
\citet{chen2012joint}.
Our results from Section~\ref{sec:truncSfs}
improve this to $O(n^2V)$.

To use
\eqref{eq:chen:sfs:formula},
we must also compute the terms
$\P(\mathbf{x} \mid \killedSize_0^\popVertex, \numAnc_0^\popVertex)
\P(\numAnc_0^\popVertex)$,
for which \citet{chen2012joint}
constructs a dynamic program,
starting at the leaf vertices
and moving up the graph.
This dynamic program essentially consists
of setting up a
Bayesian graphical model
with random variables
$\numAnc_0^\popVertex, \killedSize_0^\popVertex$
and performing belief
propagation,
which can be done via
the sum-product algorithm (``tree-peeling'')
if the population graph is a tree
\citep{pearl1982reverend, felsenstein1981evolutionary},
or via a junction tree algorithm
if not \citep{lauritzen1988local}.

The time complexity
of the algorithm
thus
depends on the topological structure
of the population graph.
In the special case where
the demographic history is a binary
tree, the tree-peeling
algorithm
computes
the values
$\P(\mathbf{x} \mid \killedSize_0^\popVertex, \numAnc_0^\popVertex)
\P(\numAnc_0^\popVertex)$
in $O(n^4V)$ time,
since the vertex $v$ has $O(n^2)$
possible states $(\killedSize_0^\popVertex, \numAnc_0^\popVertex)$,
so summing over the transitions between every pair of states costs $O(n^4)$.
Note that \citet{chen2012joint}
mistakenly states that
the computation takes $O(n^3V)$ time.
In the further special case
that the population sizes are piecewise constant,
speedups from  \citet{bryant2012inferring}
can improve this to $O(n^2 \log(n) V)$.
More specifically, \citet{bryant2012inferring}
computes the terms
$\P(\mathbf{x} \mid \killedSize_0^\popVertex, \numAnc_0^\popVertex)
\P(\numAnc_0^\popVertex)$
in $O(n^2\log(n)V)$ time
for a model with recurrent mutation,
but the results can be applied straightforwardly
here by setting the mutation rate to $0$,
thus disallowing recurrent mutation.

To summarize,
let $W$ be the time it takes to compute
\eqref{eq:chen:sfs:formula}
after the terms $\usfs_\numAnc^\popVertex(\killedSize)$
have been precomputed,
and let $\numLoci$ be the number of distinct entries
$\mathbf{x}$
for which we wish to compute 
$\usfs(\mathbf{x})$.
Then our results from Section~\ref{sec:truncSfs}
improve the computational complexity
from $O(n^5V + W \numLoci)$ to
$O(n^2V + W \numLoci)$.
In the case of a binary tree
the original algorithm of \citet{chen2012joint}
gives $W = O(n^4 V)$,
but adapting results from \citet{bryant2012inferring} improves the this to $W = O(n^2 \log(n) V)$ when the population sizes are piecewise constant.
In the following section, we introduce a new approach 
that further improves the runtime to
$W = O(n \log(n) V)$
and generalizes from
piecewise constant to arbitrary population size
histories.

\subsection{A Moran-based dynamic program\label{sec:moran}}

We describe a modified version
of the dynamic programs
from \citet{chen2012joint, bryant2012inferring}
that improves the computational
complexity of computing
$\usfs(\mathbf{x})$ for tree-shaped demographies.
The main idea is to replace
the backwards-in-time coalescent
with a forwards-in-time Moran model.

We assume the $\mathcal{D}$ populations at the present
are related by a binary rooted tree
with $\mathcal{D}$ leaves,
where each leaf represents a population
at the present,
and at each internal vertex,
a parent population
splits into
two child populations.
(Note that a non-binary tree can be represented
as a binary tree, with additional vertices of
height $0$).

Instead of working with the
multi-population coalescent directly,
we will consider a multi-population Moran
model, in which the coalescent is embedded \citep{moran1958random}.
In particular, let $\mathfrak{L}(\popVertex)$
denote the leaf populations
descended from the population $\popVertex$,
and let $n_{\popVertex} = \sum_{i \in \mathfrak{L}(\popVertex)} n_i$
be the number of present-day alleles
with ancestry in $\popVertex$.
For each population $\popVertex$ (except the root),
we construct a Moran model going \emph{forward} in time,
i.e. starting at $\truncTime_\popVertex$ and ending at $0$.
The Moran model consists of
$n_\popVertex$ lineages,
each with either an ancestral or derived allele.
Going forward in time,
every lineage copies itself onto every other
lineage at rate $\frac12 \alpha_\popVertex(\timeVar)$.
Thus, the total rate of copying events is
${n_\popVertex \choose 2} \alpha_\popVertex(\timeVar)$.
Let $\mu^\popVertex_\timeVar$ denote the number of
derived alleles at time $\timeVar$ in population $\popVertex$.
Then the transition rate of
$\mu^\popVertex_\timeVar$ when
$\mu^\popVertex_\timeVar=x$
is
$\lambda_{x \to x+1}(\timeVar) = \lambda_{x\to x-1}(\timeVar)
= \frac{x (n_\popVertex - x)}2 \alpha_\popVertex(\timeVar)$,
since there are $x (n_\popVertex - x)$
pairs of lineages with different alleles.

The coalescent is embedded within
the Moran model,
because if we trace the ancestry of genetic
material backwards in time in the Moran model,
we obtain a genealogy with the same
distribution as under the coalescent
(Theorem 1.30 of
\citet{durrett2008probability}).
Thus, we can obtain the expected
number of mutations
with derived counts $\mathbf{x}$,
by summing over which population $\popVertex$
the mutation occurred in:
\begin{eqnarray}
  \usfs(\mathbf{x})
&=&
\sum_{\popVertex}
\sum_{\killedSize=1}^{n_\popVertex}
\usfs^\popVertex_{n_\popVertex}(\killedSize)
\P(\mathbf{x} \mid \mu^\popVertex_0 = \killedSize).
\label{eq:multiSfs:sum}
\end{eqnarray}

Let $\mathbf{x}_{\popVertex} = \{x_i : i \in
\mathfrak{L}(\popVertex)\}$
denote the subsample of derived allele counts
in the populations descended from $\popVertex$.
Similarly, let $\mathbf{x}_{\popVertex}^c = \{x_i : i \notin
\mathfrak{L}(\popVertex)\}$.
Then for $\killedSize \geq 1$,
\begin{eqnarray}
  \P(\mathbf{x} \mid \mu^\popVertex_0 = \killedSize)
&=&
\begin{cases}
    \P(\mathbf{x}_{\popVertex} \mid \mu^\popVertex_0 = \killedSize),
& \text{if } \mathbf{x}_{\popVertex}^c = \mathbf{0}, \\
0,
& \text{if }\mathbf{x}_{\popVertex}^c \neq \mathbf{0}.
\end{cases}
\label{eq:bottom:likelihood}
\end{eqnarray}
So it suffices to compute
$\P(\mathbf{x}_{\popVertex} \mid \mu^\popVertex_0 = \killedSize)$
for all $\popVertex$ and $\killedSize$.
If $\popVertex$ is the $i$th leaf population,
then
$\P(\mathbf{x}_\popVertex \mid \mu^{\popVertex}_0 = \killedSize)
= \I_{\killedSize = x_i}$.
On the other hand, if $\popVertex$ is an interior vertex
with children $\popVertex_1$ and $\popVertex_2$, then
\begin{equation}
\P(\mathbf{x}_{\popVertex} \mid \mu^{\popVertex}_0 =\killedSize)
=
\sum_{\killedSize_1=0}^{n_{\popVertex_1}}
\frac{{n_{\popVertex_1} \choose \killedSize_1}
{n_{\popVertex_2} \choose \killedSize - \killedSize_1}}{{n_\popVertex \choose \killedSize} }
\P(\mathbf{x}_{\popVertex_1} \mid
\mu^{\popVertex_1}_{\truncTime_{\popVertex_1}} = \killedSize_1)
\P(\mathbf{x}_{\popVertex_2} \mid
\mu^{\popVertex_2}_{\truncTime_{\popVertex_2}} = \killedSize -
\killedSize_1), 
\label{sum:product:interior}
\end{equation}
where $\P(\mathbf{x}_{\popVertex_i} \mid
\mu^{\popVertex_i}_{\truncTime_{\popVertex_i}})$
can be computed from
\begin{eqnarray}
    \P(\mathbf{x}_{\popVertex} \mid \mu^{\popVertex}_{\truncTime_\popVertex} =
  \killedSize)
&=&
\sum_{j=0}^{n_{\popVertex}}
\P(\mathbf{x}_\popVertex \mid \mu^{\popVertex}_{0} = j)
\P(\mu^{\popVertex}_{0} =
  j \mid \mu^{\popVertex}_{\truncTime_\popVertex} =
  \killedSize).
\label{eq:sfs:transition:1}
\end{eqnarray}
To compute
the transition probability
$\P(\mu^{\popVertex}_{0} =
  j \mid \mu^{\popVertex}_{\truncTime_\popVertex} =
  \killedSize)$,
note that the transition rate matrix
of $\mu^\popVertex_\timeVar$
can be written as $Q^{(\popVertex)} \alpha(\timeVar)$,
where $Q^{(\popVertex)} = (q^{(\popVertex)}_{ij})_{0\leq i,j
  \leq n_{\popVertex}}$ is a $(n+1)\times(n+1)$ matrix
with
\begin{eqnarray*}
  q^{(\popVertex)}_{ij}
&=&
\begin{cases}
  - i (n_\popVertex - i), & \text {if } i=j, \\
  \frac{1}2 i (n_\popVertex-i), & \text{if } |j-i|=1, \\
  0, & \text{else,}
\end{cases}
\end{eqnarray*}
so then the transition probability is given by the
matrix exponential
\begin{eqnarray}
  \P(\mu^{\popVertex}_{0} =
  j \mid \mu^{\popVertex}_{\truncTime_\popVertex} =
  \killedSize)
&=&
( e^{Q^{(\popVertex)} \int_0^{\truncTime_\popVertex}
  \alpha_\popVertex(\timeVar) d\timeVar} )_{\killedSize,j} .
\label{eq:sfs:transition:exp}
\end{eqnarray}

Thus, the joint SFS $\usfs(\mathbf{x})$
can be computed using
\eqref{eq:multiSfs:sum}
and \eqref{eq:bottom:likelihood},
with
$\P(\mathbf{x}_{\popVertex} \mid \mu^\popVertex_0 = \killedSize)$
given by recursively
computing
\eqref{sum:product:interior},
\eqref{eq:sfs:transition:1},
and \eqref{eq:sfs:transition:exp},
in a depth-first search
on the population tree
(i.e.
Felsenstein's tree-peeling algorithm,
or the sum-product algorithm 
for belief propagation).

We now consider the computational
complexity associated with each vertex
$\popVertex$.
Equations
\eqref{sum:product:interior}
and
\eqref{eq:sfs:transition:1}
each have $O(n_\popVertex)$
terms,
and must be solved for $O(n_\popVertex)$
values of $\killedSize$;
so naively, each vertex costs
$O(n_\popVertex^2)$ time.
However, we can improve \eqref{sum:product:interior}
to $O(n_\popVertex \log(n_\popVertex))$
and \eqref{eq:sfs:transition:1}
to $O(n_\popVertex)$,
using
essentially the same speedups as in \citet{bryant2012inferring}.
Letting
$\tilde\likelihood^\popVertex_\timeVar(k) = {n_{\popVertex} \choose k}
\P(\mathbf{x}_\popVertex \mid \mu^{\popVertex}_{\timeVar} = k)$,
\eqref{sum:product:interior} can be written as a convolution
\begin{eqnarray}
  \tilde\likelihood^\popVertex_0 =
  \tilde\likelihood^{\popVertex_1}_{\truncTime_{\popVertex_1}} *
  \tilde\likelihood^{\popVertex_2}_{\truncTime_{\popVertex_2}},
\label{eq:fftconvolve}
\end{eqnarray}
which can be computed in $O(n_\popVertex \log(n_\popVertex))$ time via the
fast Fourier transform \citep{cooley1965algorithm},
since $\mathcal{F} \likelihood^\popVertex_0
= \left(\mathcal{F}
  \likelihood^{\popVertex_1}_{\truncTime_{\popVertex_1}}\right)
\left(\mathcal{F}
  \likelihood^{\popVertex_2}_{\truncTime_{\popVertex_2}}\right)$,
where $\mathcal{F}$ is the discrete Fourier transform.
Similarly, letting $\likelihood^\popVertex_\timeVar(k) = \tilde\likelihood^\popVertex_\timeVar(k)/ {n_{\popVertex} \choose k}$,
\eqref{eq:sfs:transition:1} turns into
\begin{eqnarray}
\likelihood^\popVertex_{\truncTime_\popVertex} &=&
e^{\left({Q^{(\popVertex)} \int_0^{\truncTime_\popVertex}
  \alpha_\popVertex(\timeVar) d\timeVar}\right)} 
\likelihood^\popVertex_{0},
\label{eq:exp:action}
\end{eqnarray}
and this costs $O(n_\popVertex)$
by the sparsity of 
$Q^{(\popVertex)}$,
using results 
for computing the action of sparse matrix exponentials
\citep{sidje1998expokit, al-mohy2011computing}.  Transforming between $\tilde\likelihood^\popVertex_{\truncTime_\popVertex}$ and $\likelihood^\popVertex_{\truncTime_\popVertex}$ takes $O(n_\popVertex)$ time.

The computational complexity associated with a single vertex
$\popVertex$ is thus $O(n_\popVertex \log(n_\popVertex))$. Therefore,
computing the joint SFS entry $\usfs(\mathbf{x})$ for $\numLoci$ distinct values of $\mathbf{x}$ takes
$O(n^2V + n \log(n) V \numLoci)$ time for a binary population tree
with arbitrary population size functions and no migration. This is a
substantial improvement over the $O(n^5V + n^4V \numLoci)$ complexity of
\citet{chen2012joint}, and the $O(n^2 \log(n) V \numLoci)$ complexity of
\citet{bryant2012inferring}.  Similar to \citet{chen2012joint}, our approach
has the benefit of easily generalizing to arbitrary population size
histories, not just piecewise constant sizes.

\section{Results}\label{sec:results}

\begin{figure}[t]
  \centering
  \includegraphics[width=0.8\textwidth]{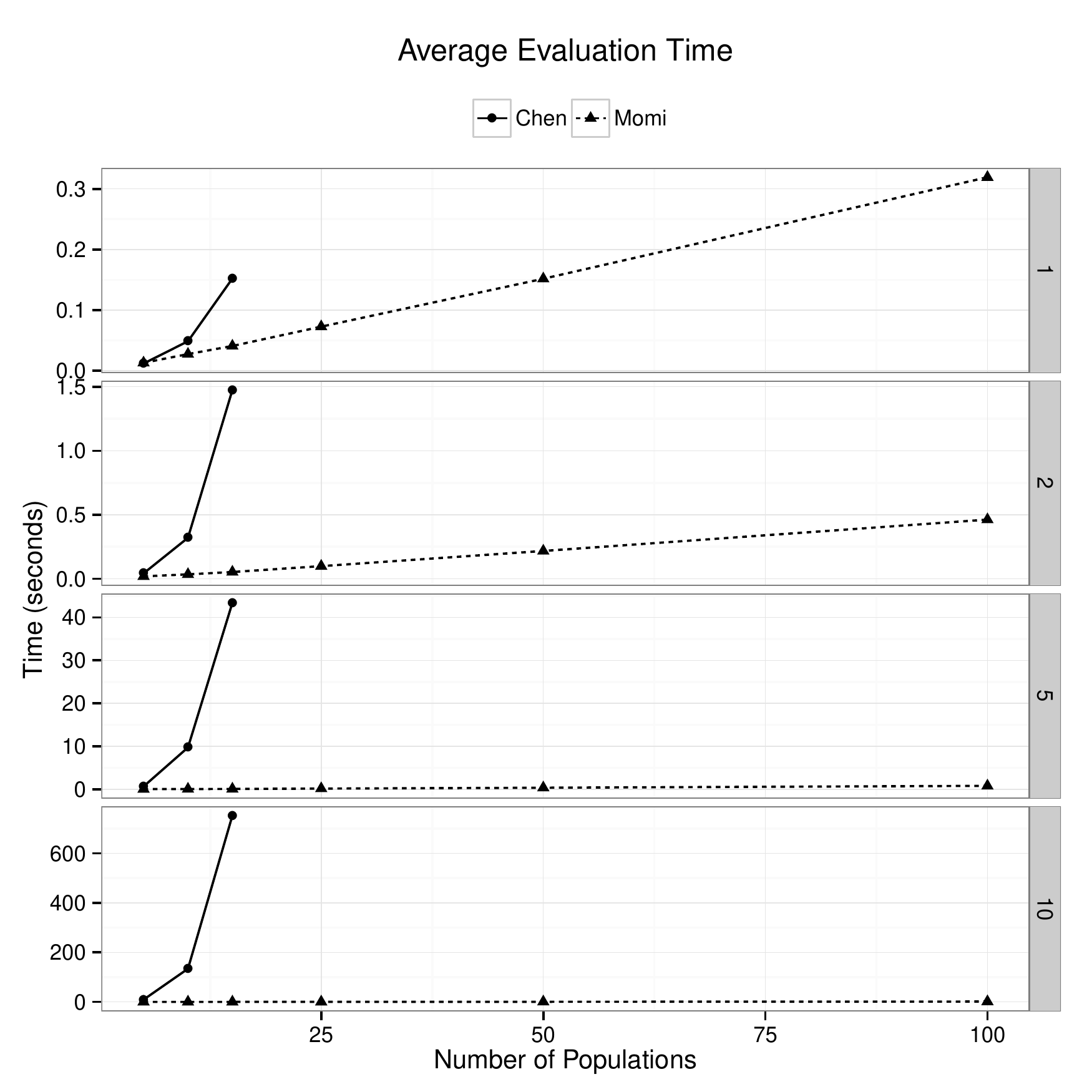}
  \caption{Average computation time per joint SFS entry.
      For each combination of the number $\mathcal{D}$ of populations and the sample size $n / \mathcal{D}$ per population, 
      we generated 20 random datasets, each under a demographic history that is a random binary tree.
      The expected joint SFS for the resulting segregating sites were then computed using our method (\emph{momi}) and that of \citet{chen2012joint}. 
      Average runtime (in seconds) per joint SFS entry  is plotted on the $y$-axis, with each panel
      corresponding to a different value of $n / \mathcal{D}$. 
      As the plots show, our algorithm is orders of magnitude faster than Chen's.
       Due to its significantly
      increased runtime, we were able to run Chen's method only up to $\mathcal{D}=15$.
       }
  \label{fig:timing}
\end{figure}

\begin{figure}[t]
  \centering
  \includegraphics[width=0.7\textwidth]{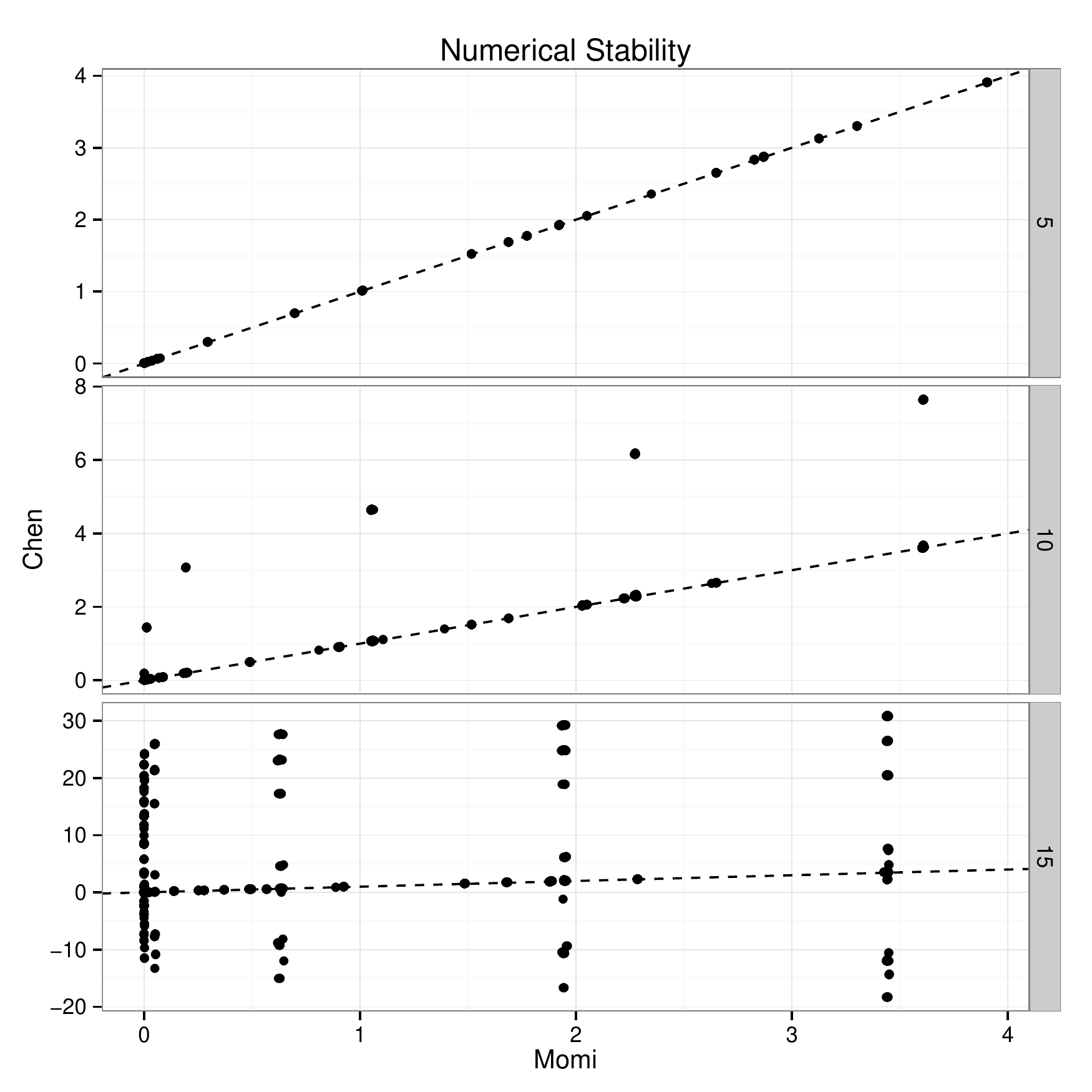}
  \caption{Numerical stability of the two algorithms.
      The plot compares the numerical values returned by our method (\emph{momi}) and
  Chen's method, for the simulations described in Figure \ref{fig:timing}.
  The three panels on the $y$-axis correspond to $\mathcal{D}\in \{5,10,15\}$. 
  To adequately illustrate the observed range of numerical values, the SFS values were transformed via the map $z \mapsto {\rm sign}(z) \log_{10}(1 + |z|)$; the
  dashed line represents the identity $y=x$.  The two methods agree for
$\mathcal{D} \le 5$, but Chen's method displays considerable numerical instability for
$\mathcal{D} \geq 10$.}
\label{fig:accuracy}
\end{figure}

We implemented our formulas and algorithm in Python, using
the Python packages \emph{numpy} and \emph{scipy}. We also implemented the formulas from \citet{chen2012joint,
  chen2013intercoalescence}, and compared the performance of the
two algorithms on simulated data.

We simulated data for demographic trees with $\mathcal{D}
\in \{5,10,15,25,50,100\}$ populations at the present, and
$\frac{n}{\mathcal{D}} \in \{1, 2, 5, 10\}$ individuals per population.
For each value of $n, \mathcal{D}$, we used the program \emph{scrm}
\citep{Staab2015} to generate 20 random datasets, 
each with a demographic history that is a random binary tree.

In Figure~\ref{fig:timing}, we
compare the running time of the original algorithm of
\citet{chen2012joint, chen2013intercoalescence} against 
our new algorithm that utilizes the formulas
for $\usfs_n^\truncTime(\killedSize)$ presented in Section~\ref{sec:truncSfs}
and our new Moran-based approach described in Section~\ref{sec:moran}.
We find our algorithm to be orders of magnitude faster; the difference
is especially pronounced as the number $\mathcal{D}$ of populations grows.
Note that, due to the increased running time, we were only able to run
Chen's algorithm to completion for $\mathcal{D} \le 15$.

In Figure~\ref{fig:accuracy}, we compare the accuracy of the two
algorithms. The figure compares the SFS entries returned by the
two methods across a subset of the simulations depicted in Figure~\ref{fig:timing}.
To adequately capture the large range of numerical values returned by the Chen
method, we transformed each SFS entry using the transformation $z \mapsto {\rm sign}(z) \log_{10}(1 + |z|)$.
The line $y=x$ is also plotted; points falling on the line depict the SFS
entries where both methods agreed. All negative return values represent numerical errors.  The two methods agree for
$\mathcal{D} \le 5$, but Chen's algorithm displays considerable numerical instability for
$\mathcal{D} = 10$ and higher.

\section{Proofs}\label{sec:proofs}
In this section, we provide proofs of the mathematical results presented in earlier sections.

\subsection{A recursion for efficiently computing $\P_\nu(\AncCoal_\truncTime = \numAnc)$}
\label{sec:ancProbs}

We describe how to compute
$\P_\nu(\AncCoal_\truncTime = \numAnc)$, for all
values of $\numAnc \leq \nu \leq n$, in
$O(n^2)$ time.
First, note that
\begin{eqnarray*}
\lefteqn{\P_{\nu-1}(\AncCoal_\truncTime = \numAnc)}\\
&=&
\P_{\nu}(\AncCoal_\truncTime = \numAnc+1, \{\nu\} \in
\Coal_\truncTime)
+
\P_{\nu}(\AncCoal_\truncTime = \numAnc,
\{\nu\} \notin \Coal_\truncTime)
\\
&=&
\frac{(\numAnc+1) \polyaProb{\nu}{\numAnc+1}{1}{1}}{{\nu \choose 1}}
\P_{\nu}(\AncCoal_\truncTime = \numAnc+1)
+
\left(1 - 
\frac{\numAnc \polyaProb{\nu}{\numAnc}{1}{1}}{{\nu \choose 1}}
\right)
\P_{\nu}(\AncCoal_\truncTime = \numAnc)
\\
&=&
\frac{(\numAnc+1)(\numAnc)}{\nu(\nu-1)}
\P_{\nu}(\AncCoal_\truncTime = \numAnc+1)
+
\left(1 - \frac{\numAnc(\numAnc-1)}{\nu(\nu-1)}\right)
\P_{\nu}(\AncCoal_\truncTime = \numAnc).
\end{eqnarray*}
Rearranging, we get the recursion
\begin{equation}
  \P_\nu(\AncCoal_\truncTime = \numAnc)
=
\frac{1}{1 - \frac{\numAnc(\numAnc-1)}{\nu(\nu-1)}}
\left[
\P_{\nu-1}(\AncCoal_\truncTime = \numAnc)
-
\frac{(\numAnc+1)(\numAnc)}{\nu(\nu-1)}
\P_{\nu}(\AncCoal_\truncTime = \numAnc+1)
\right]
\end{equation}
with base cases
\begin{eqnarray*}
  \P_\nu(\AncCoal_\truncTime = \nu)
&=&
e^{-{\nu \choose 2} \int_0^\truncTime \alpha(\timeVar) d \timeVar}.
\end{eqnarray*}
So after solving $\int_0^\truncTime \alpha(\timeVar) d \timeVar$,
we can use the recursion and memoization to
solve for
all of the $O(n^2)$ terms
$\P_\nu(\AncCoal_\truncTime = \numAnc)$
in $O(n^2)$ time.
In particular, in the case of
constant population size, $\alpha(\timeVar) = \alpha$,
the base case is given by
\begin{eqnarray*}
  \P_\nu(\AncCoal_\truncTime = \nu)
&=&
e^{-{\nu \choose 2} \alpha \truncTime},
\end{eqnarray*}
and in the case of an exponentially growing population size,
$\alpha(\timeVar) = \alpha(\truncTime) e^{\beta (\truncTime -
  \timeVar)}$,
the base case is given by
\begin{eqnarray*}
  \P_\nu(\AncCoal_\truncTime = \nu)
&=&
e^{-{\nu \choose 2}\alpha(\truncTime) (e^{\beta \truncTime} - \frac1\beta)}.
\end{eqnarray*}

\subsection{Proof of \lemref{thm:tsfs:tmrca}}
\label{proof:tsfs:tmrca}

Let $\tmrca$ denote the time to the most recent common ancestor of the
sample. We first note that
\begin{eqnarray*}
  \usfs_n^\truncTime(n) &=& \truncTime - \E_n[\tmrca \wedge \truncTime],
\end{eqnarray*}
since the branch length subtending the whole sample is the time between $\truncTime$ and $\tmrca$.

Next, note that $\frac{\theta}2 \E_n[\tmrca \wedge \truncTime]$ is equal to 
the number of polymorphic mutations in $[0,\truncTime)$ where the individual ``1'' is derived.
This is because, as we trace the ancestry of ``1'' backwards in time,
all mutations hitting the lineage below $\tmrca$ are polymorphic,
while all mutations hitting above $\tmrca$ are monomorphic.

The expected number of polymorphic mutations with ``1'' derived
is also equal to $\frac{\theta}2 \sum_{\killedSize=1}^{n-1}
\frac{\killedSize}n \usfs_n^\truncTime(\killedSize)$, since if a
mutation has $\killedSize$ derived leaves, the chance that ``1''
is in the derived set is $\frac{\killedSize}n$. Thus,
\begin{eqnarray*}
  \E_n[\tmrca \wedge \truncTime]
&=&
\sum_{\killedSize=1}^{n-1} \frac{\killedSize}n \usfs_n^\truncTime(\killedSize),
\end{eqnarray*}
which completes the proof.

\subsection{Proof of \lemref{thm:tsfs:recurrence}} \label{proof:tsfs:recurrence}
We first note that
\begin{eqnarray*}
\lefteqn{\P_n(\mutSet^\truncTime = \{1,\ldots,\killedSize\})}\\
&=&
  \P_{n+1}(\mutSet^\truncTime = \{1,\ldots,\killedSize\})
+
  \P_{n+1}(\mutSet^\truncTime = \{1,\ldots,\killedSize, n+1\}).
\end{eqnarray*}
By exchangeability, we have $\P_n(\mutSet^\truncTime = \killedSet) =
\frac{\theta}2  \frac{\usfs_n^\truncTime(|\killedSet|)}{{n\choose |\killedSet|}}
+ o(\theta)$ for all $\killedSet \subseteq \{1,\ldots,n\}$, so
\begin{eqnarray*}
  \frac1{{n \choose \killedSize}}
\usfs_n^\truncTime(\killedSize)
&=&
\frac1{{n+1 \choose \killedSize}} \usfs_{n+1}^\truncTime(\killedSize)
+
\frac1{{n+1 \choose \killedSize+1}} \usfs_{n+1}^\truncTime(\killedSize+1).
\end{eqnarray*}
Multiplying both sides by ${n \choose \killedSize}$ gives
\begin{eqnarray*}
\usfs_n^\truncTime(\killedSize)
&=&
\frac{n-k+1}{n+1} \usfs_{n+1}^\truncTime(\killedSize)
+
\frac{k+1}{n+1} \usfs_{n+1}^\truncTime(\killedSize+1).
\end{eqnarray*}
  
\subsection{Proof of \lemref{thm:truncated:polanski:kimmel}\label{proof:truncated:polanski:kimmel}}

Let $\alpha^*(\timeVar)$ denote the inverse population size history given by
\begin{eqnarray*}
  \alpha^*(\timeVar) &=&
  \begin{cases}
    \alpha(\timeVar) & \text{if } \timeVar < \truncTime \\
    \infty & \text{if } \timeVar \geq \truncTime.
  \end{cases}
\end{eqnarray*}
So the demographic history with population size $\frac1{\alpha^*(\timeVar)}$
agrees with the original history up to time $\truncTime$, at which
point the population size drops to $0$, and all lineages instantly coalesce
into a single lineage with probability $1$.

Let $T_{m,*}$ denote the amount of time there are $m$ ancestral lineages for the coalescent with size history $\frac1{\alpha^*(\timeVar)}$. Similarly, let 
$\usfs_{n,*}(\killedSize)$ denote the SFS under the size history $\frac1{\alpha^*(\timeVar)}$. Then from the result of \citet{polanski2003new},
\begin{equation*}
  \usfs_{n,*}(\killedSize) = \sum_{m=2}^n W_{n,k,m} \E_m[T_{m,*}].
\end{equation*}
Note that for $m > 1$, we almost surely have $T_{m,*} = T^\truncTime_{m,*}$,
i.e. the intercoalescence time equals its truncated version,
since all lineages coalesce instantly at $\truncTime$ with probability $1$.
Thus, $\E_m[T_{m,*}] = \E_m[T^\truncTime_{m,*}]$.
Similarly, for $\killedSize < n$, $\usfs_{n,*}(\killedSize) = \usfs^\truncTime_{n,*}(\killedSize)$, i.e. the SFS equals the truncated SFS, because the probability of a polymorphic mutation occurring in $[\tau, \infty)$ is $0$.

Finally, note that $\E_m[T^\truncTime_{m,*}] = \E_m[T^\truncTime_m]$ and $\usfs^\truncTime_{n,*}(\killedSize) = \usfs^\truncTime_{n}(\killedSize)$, because $\alpha(\timeVar)$ and $\alpha^*(\timeVar)$ are identical on $[0,\truncTime)$.

\subsection{Proof of Proposition~\ref{thm:csfs:ancKill}\label{proof:csfs:ancKill}}

We start by showing that
$    \P_n(\AncKill_\truncTime =
    \numAnc) = \P_n(\AncCoal_\truncTime = \numAnc) + O(\theta)$.
Let $\truncWait{i}(\CKill) = \int_0^\truncTime \I_{\AncKill_\timeVar =
  i} d \timeVar$ 
denote the amount of time where $\CKill$ has $i$ unkilled lineages.
Let $p$ denote the probability density function.
For $(\timeVar_n,\ldots, \timeVar_\numAnc)$ with $\sum \timeVar_i = \truncTime$,
we have
\begin{eqnarray*}
&&
  p(\truncWait{n}(\CKill) = \timeVar_n, \ldots, \truncWait{\numAnc}(\CKill)
  = \timeVar_{\numAnc})
\\
&=&
e^{-\lambda^\CKill_{\numAnc,\numAnc-1} \timeVar_\numAnc}
\prod_{i=\numAnc+1}^n \lambda^\CKill_{i,i-1} e^{-\lambda^\CKill_{i,i-1}
  \timeVar_i}
\\
&=&
e^{-\left({\numAnc \choose 2} \alpha + \frac{\numAnc \theta}2 \right) \timeVar_\numAnc}
\prod_{i=\numAnc+1}^n \left({i \choose 2} \alpha + \frac{i \theta}2\right) e^{-\left({i \choose 2} \alpha + \frac{i \theta}2\right)
  \timeVar_i}
\\
&=&
e^{-{\numAnc \choose 2} \alpha \timeVar_\numAnc}
\prod_{i=\numAnc+1}^n {i \choose 2} \alpha e^{-{i \choose 2} \alpha
  \timeVar_i}
+ O(\theta)
\\
&=&
p(\truncWait{n} = \timeVar_n, \ldots, \truncWait{\numAnc} = \timeVar_\numAnc)
+ O(\theta),
\end{eqnarray*}
and so
\begin{eqnarray*}
  \lim_{\theta \to 0} \P_n(\AncKill_\truncTime = \numAnc)
&=&
\lim_{\theta \to 0} \int_{\sum \timeVar_i = \truncTime}
  p(\truncWait{n}(\CKill) = \timeVar_n, \ldots, \truncWait{\numAnc}(\CKill)
  = \timeVar_{\numAnc}) d \mathbf{\timeVar}
\\
&=&
\int_{\sum \timeVar_i = \truncTime}
p(\truncWait{n} = \timeVar_n, \ldots, \truncWait{\numAnc} = \timeVar_\numAnc)
   d \mathbf{\timeVar}
\\
&=&
\P_n(\AncCoal_\truncTime = \numAnc).
\end{eqnarray*}
where we can exchange the limit
and the integral by the Bounded Convergence
Theorem,
because
$  p(\truncWait{n}(\CKill) = \timeVar_n, \ldots, \truncWait{\numAnc}(\CKill)
  = \timeVar_{\numAnc}) \leq \prod_{i=\numAnc+1}^n \left({i \choose 2} \alpha
    + \frac{i}2 \right)$ for $\theta \leq 1$.

Thus we have
\begin{eqnarray*}
      \P_n(|\mutSet^\truncTime| = \killedSize, \AncKill_\truncTime =
    \numAnc)
&=&
    \P_n(|\mutSet^\truncTime| = \killedSize \mid \AncKill_\truncTime =
    \numAnc)
    \P_n(\AncKill_\truncTime =
    \numAnc)
\\
&=&
\left( \frac{\theta}2   
\usfs_n^\truncTime(\killedSize \mid \AncKill_\truncTime = \numAnc)
+ o(\theta) \right)
\left( \P_n(\AncCoal_\truncTime = \numAnc) + O(\theta) \right)
\\
&=&
\frac{\theta}2   
\usfs_n^\truncTime(\killedSize \mid \AncKill_\truncTime = \numAnc)
\P_n(\AncCoal_\truncTime = \numAnc) + o(\theta),
\end{eqnarray*}
which proves the first part of the proposition.

We next solve for
$\usfs_n^\truncTime(\killedSize
  \mid \AncKill_\truncTime = \numAnc)$,
the first
order Taylor series
coefficient for
$\P_n(|\mutSet^\truncTime| = \killedSize \mid \AncKill_\truncTime =
\numAnc)$
in the mutation rate $\frac{\theta}2$.

When there are $i$ unkilled lineages,
the probability that the next event
is a killing event is
$\frac{\theta}{\alpha(i-1) + \theta} = \frac{\theta}{\alpha(i-1)} +
o(\theta)$.
Given that the event is a killing,
the chance that the killed lineage
has $\killedSize$ leaf descendants
is $\polyaProb{n}{i}{\killedSize}{1}$.
So summing over $i$,
and dividing out the mutation rate $\frac{\theta}2$,
we get
  \begin{eqnarray*}
    \usfs_n^\truncTime(\killedSize \mid \AncKill_\truncTime = \numAnc)
    &=&
    \frac2\alpha \sum_{i=m+1}^{n-\killedSize+1} \frac1{i-1}
    \polyaProb{n}{i}{\killedSize}{1}
\\
&=&
\frac2\alpha
\sum_{i=m+1}^{n-\killedSize+1} \frac1{i-1} \frac{{n - \killedSize - 1 \choose i -
    2}}{{n - 1 \choose i -1}}
\\
&=&
\frac2\alpha
\sum_{i=m+1}^{n-\killedSize+1}
\frac1{i-1}
\frac{(n-k-1)! (i-1)! (n-i)!}{ (i-2)! (n-k-i+1)! (n-1)!}
\\
&=&
\frac{2 (n-k-1)!}{\alpha(n-1)!}
\sum_{i=m+1}^{n-\killedSize+1}
\frac{ (n-i)!}{  (n-k-i+1)! }
\\
&=&
\frac{2 (n-k-1)!}{\alpha (n-1)!}
\sum_{j=0}^{n-\killedSize-\numAnc}
\frac{ (j+\killedSize-1)!}{  j! }
\\
&=&
\frac{2}{\alpha \killedSize {n - 1 \choose \killedSize}}
\sum_{j=0}^{n-\killedSize-\numAnc}
{j + \killedSize - 1 \choose j}
\\
&=&
\frac{2}{\alpha \killedSize}
\frac{{n-\numAnc \choose \killedSize}}{ {n - 1 \choose \killedSize}},
  \end{eqnarray*}
where we made the change of variables
$j = n - \killedSize - i + 1$,
and where the final line
follows from repeated
application of the combinatorial identity
${a \choose b} = {a-1 \choose b} + {a-1 \choose b-1}$.

\subsubsection{Alternative proof for $\usfs_n^\truncTime(\killedSize
  \mid \AncKill_\truncTime = \numAnc)$ via the Chinese Restaurant Process}

We sketch an alternative proof of the expression for
$\usfs_n^\truncTime(\killedSize
  \mid \AncKill_\truncTime = \numAnc)$,
using the Chinese Restaurant Process.

Consider the coalescent with killing going forward in time
(towards the present),
and only looking at it when
the number of individuals increases.
Then when there are $i$ lineages,
a new mutation occurs with probability
$\frac{\theta}{\alpha i + \theta} = \frac{\theta/\alpha}{i +
  \theta/\alpha}$,
and each lineage branches with probability
$\frac{\alpha }{\alpha i + \theta} = \frac{1}{i + \theta/\alpha}$.
Thus, conditional on $\AncKill_\truncTime = \numAnc$,
the distribution on $\CKill_\truncTime$ is given by a Chinese Restaurant Process
\citep{aldous1985exchangeability},
starting with $\numAnc$ tables each with $1$ person, and with new tables founded
with parameter $\theta / \alpha$.

Let 
$\rise{x}{i} = x (x+1) \cdots (x+i-1)$
denote the rising factorial.
If there is a single mutation with
$\killedSize$ descendants,
then there are ${n - \numAnc \choose \killedSize}$
ways to pick which of the $n-\numAnc$
events involve mutant lineages.
The probability of a particular such ordering
is
  \begin{eqnarray*}
    \frac{\theta}\alpha \frac{\rise{1}{\killedSize} \rise{\numAnc}{n - \killedSize
        - \numAnc}}{\rise{\numAnc+\theta/\alpha}{n-\numAnc}}
    &=& \frac{\theta}\alpha \frac{(\killedSize -1)! (n - \killedSize - 1)! /
      \numAnc!}{(n-1)! / \numAnc!} + o(\theta).
  \end{eqnarray*}
  Summing over all ${n - \numAnc \choose \killedSize}$ orderings,
  and dividing by $\frac{\theta}2$, yields
  \begin{align*}
    \usfs_n^\truncTime(\killedSize \mid \AncKill_\truncTime = \numAnc)
    &=
    \frac{2}\alpha {n - \numAnc \choose \killedSize} \frac{(\killedSize -1)! (n - \killedSize - 1)! / \numAnc!}{(n-1)! / \numAnc!}.
  \end{align*}

\bibliographystyle{imsart-nameyear}   
\bibliography{yss-group}

\end{document}

%% file: figures/notation.pdf_tex
\begingroup%
  \makeatletter%
  \providecommand\color[2][]{%
    \errmessage{(Inkscape) Color is used for the text in Inkscape, but the package 'color.sty' is not loaded}%
    \renewcommand\color[2][]{}%
  }%
  \providecommand\transparent[1]{%
    \errmessage{(Inkscape) Transparency is used (non-zero) for the text in Inkscape, but the package 'transparent.sty' is not loaded}%
    \renewcommand\transparent[1]{}%
  }%
  \providecommand\rotatebox[2]{#2}%
  \ifx\svgwidth\undefined%
    \setlength{\unitlength}{495.74467773bp}%
    \ifx\svgscale\undefined%
      \relax%
    \else%
      \setlength{\unitlength}{\unitlength * \real{\svgscale}}%
    \fi%
  \else%
    \setlength{\unitlength}{\svgwidth}%
  \fi%
  \global\let\svgwidth\undefined%
  \global\let\svgscale\undefined%
  \makeatother%
  \begin{picture}(1,0.6)%
    \put(0,0){\includegraphics[width=\unitlength]{figures/notation.pdf}}%
    \put(0.19479264,0.00188099){\color[rgb]{0,0,0}\makebox(0,0)[lb]{\smash{$1$}}}%
    \put(0.33873157,0.00188099){\color[rgb]{0,0,0}\makebox(0,0)[lb]{\smash{$2$}}}%
    \put(0.46837974,0.00188099){\color[rgb]{0,0,0}\makebox(0,0)[lb]{\smash{$3$}}}%
    \put(0.61446097,0.00188099){\color[rgb]{0,0,0}\makebox(0,0)[lb]{\smash{$4$}}}%
    \put(0.7598066,0.00091407){\color[rgb]{0,0,0}\makebox(0,0)[lb]{\smash{$5$}}}%
    \put(-0.03,0.095){\color[rgb]{0,0,0}\makebox(0,0)[lt]{\begin{minipage}{0.24501158\unitlength}\raggedright $t=0$\\ \end{minipage}}}%
    \put(-0.03,0.495){\color[rgb]{0,0,0}\makebox(0,0)[lt]{\begin{minipage}{0.1784124\unitlength}\raggedright $t=\truncTime$\end{minipage}}}%
    \put(1.03000587,0.125){\color[rgb]{0,0,0}\makebox(0,0)[lb]{\smash{$T_5^\truncTime(=T_5)$}}}%
    \put(1.03000587,0.275){\color[rgb]{0,0,0}\makebox(0,0)[lb]{\smash{$T_4^\truncTime(=T_4)$}}}%
    \put(1.03000587,0.41){\color[rgb]{0,0,0}\makebox(0,0)[lb]{\smash{$T_3^\truncTime(<T_3)$}}}%
    \put(0.16,0.50890163){\color[rgb]{0,0,0}\makebox(0,0)[lb]{\smash{$\Coal_\truncTime=\{\{1,2\},\{3\},\{4,5\}\}, \AncCoal_\truncTime=3$}}}%
    \put(0.35,0.58){\color[rgb]{0,0,0}\makebox(0,0)[lb]{\smash{$\mutSet^\truncTime = \{1,4,5\}$}}}%
  \end{picture}%
\endgroup%

%% file: figures/coalkill.pdf_tex
\begingroup%
  \makeatletter%
  \providecommand\color[2][]{%
    \errmessage{(Inkscape) Color is used for the text in Inkscape, but the package 'color.sty' is not loaded}%
    \renewcommand\color[2][]{}%
  }%
  \providecommand\transparent[1]{%
    \errmessage{(Inkscape) Transparency is used (non-zero) for the text in Inkscape, but the package 'transparent.sty' is not loaded}%
    \renewcommand\transparent[1]{}%
  }%
  \providecommand\rotatebox[2]{#2}%
  \ifx\svgwidth\undefined%
    \setlength{\unitlength}{495.74467773bp}%
    \ifx\svgscale\undefined%
      \relax%
    \else%
      \setlength{\unitlength}{\unitlength * \real{\svgscale}}%
    \fi%
  \else%
    \setlength{\unitlength}{\svgwidth}%
  \fi%
  \global\let\svgwidth\undefined%
  \global\let\svgscale\undefined%
  \makeatother%
  \begin{picture}(1,0.53)%
    \put(0,0){\includegraphics[width=\unitlength]{figures/coalkill.pdf}}%
    \put(0.19479264,0.001881){\color[rgb]{0,0,0}\makebox(0,0)[lb]{\smash{$1$}}}%
    \put(0.33873157,0.001881){\color[rgb]{0,0,0}\makebox(0,0)[lb]{\smash{$2$}}}%
    \put(0.46837974,0.001881){\color[rgb]{0,0,0}\makebox(0,0)[lb]{\smash{$3$}}}%
    \put(0.61446097,0.001881){\color[rgb]{0,0,0}\makebox(0,0)[lb]{\smash{$4$}}}%
    \put(0.7598066,0.00091408){\color[rgb]{0,0,0}\makebox(0,0)[lb]{\smash{$5$}}}%
    \put(-0.03,0.095){\color[rgb]{0,0,0}\makebox(0,0)[lt]{\begin{minipage}{0.24501158\unitlength}\raggedright $t=0$\\ \end{minipage}}}%
    \put(-0.03,0.495){\color[rgb]{0,0,0}\makebox(0,0)[lt]{\begin{minipage}{0.1784124\unitlength}\raggedright $t=\truncTime$\end{minipage}}}%
    \put(0.14928153,0.53){\color[rgb]{0,0,0}\makebox(0,0)[lb]{\smash{$\CKill_\truncTime = \{\underline{\{1\}}, \{2\}, \{3\}, \underline{\{4,5\}}\}$}}}%
  \end{picture}%
\endgroup%

%% file: figures/demography3.pdf_tex
\begingroup%
  \makeatletter%
  \providecommand\color[2][]{%
    \errmessage{(Inkscape) Color is used for the text in Inkscape, but the package 'color.sty' is not loaded}%
    \renewcommand\color[2][]{}%
  }%
  \providecommand\transparent[1]{%
    \errmessage{(Inkscape) Transparency is used (non-zero) for the text in Inkscape, but the package 'transparent.sty' is not loaded}%
    \renewcommand\transparent[1]{}%
  }%
  \providecommand\rotatebox[2]{#2}%
  \ifx\svgwidth\undefined%
    \setlength{\unitlength}{1269.08769531bp}%
    \ifx\svgscale\undefined%
      \relax%
    \else%
      \setlength{\unitlength}{\unitlength * \real{\svgscale}}%
    \fi%
  \else%
    \setlength{\unitlength}{\svgwidth}%
  \fi%
  \global\let\svgwidth\undefined%
  \global\let\svgscale\undefined%
  \makeatother%
  \begin{picture}(1,0.37013901)%
    \put(0,0){\includegraphics[width=\unitlength]{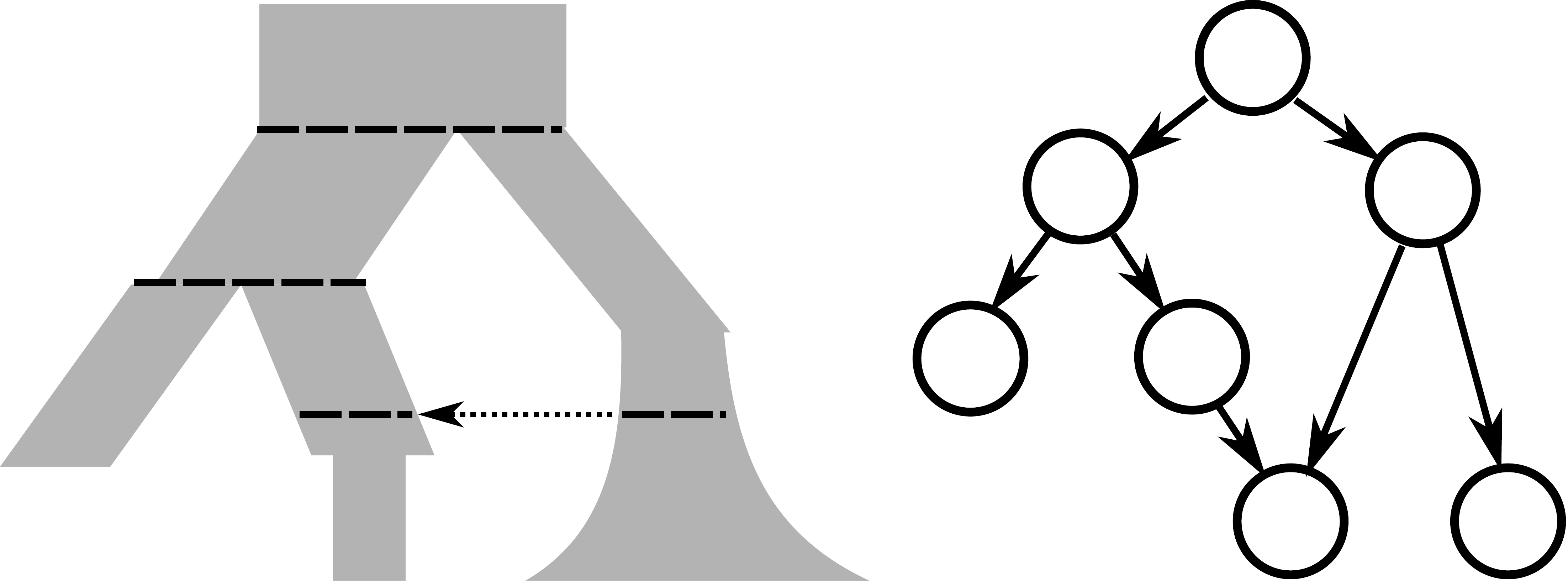}}%
    \put(0.81101636,0.02881753){\color[rgb]{0,0,0}\makebox(0,0)[lb]{\smash{$v_6$}}}%
    \put(0.94969866,0.02881753){\color[rgb]{0,0,0}\makebox(0,0)[lb]{\smash{$v_7$}}}%
    \put(0.78529959,0.32435733){\color[rgb]{0,0,0}\makebox(0,0)[lb]{\smash{$v_1$}}}%
    \put(0.67518465,0.24197884){\color[rgb]{0,0,0}\makebox(0,0)[lb]{\smash{$v_2$}}}%
    \put(0.60575432,0.13204936){\color[rgb]{0,0,0}\makebox(0,0)[lb]{\smash{$v_3$}}}%
    \put(0.74812268,0.13376438){\color[rgb]{0,0,0}\makebox(0,0)[lb]{\smash{$v_4$}}}%
    \put(0.89465549,0.23960687){\color[rgb]{0,0,0}\makebox(0,0)[lb]{\smash{$v_5$}}}%
    \put(0.24580254,0.32435733){\color[rgb]{0,0,0}\makebox(0,0)[lb]{\smash{$v_1$}}}%
    \put(0.16281454,0.22375477){\color[rgb]{0,0,0}\makebox(0,0)[lb]{\smash{$v_2$}}}%
    \put(0.06625724,0.13204936){\color[rgb]{0,0,0}\makebox(0,0)[lb]{\smash{$v_3$}}}%
    \put(0.20587697,0.12814781){\color[rgb]{0,0,0}\makebox(0,0)[lb]{\smash{$v_4$}}}%
    \put(0.394361,0.17862517){\color[rgb]{0,0,0}\makebox(0,0)[lb]{\smash{$v_5$}}}%
    \put(0.21988352,0.03607726){\color[rgb]{0,0,0}\makebox(0,0)[lb]{\smash{$v_6$}}}%
    \put(0.4234606,0.03752919){\color[rgb]{0,0,0}\makebox(0,0)[lb]{\smash{$v_7$}}}%
  \end{picture}%
\endgroup%